\documentclass[11pt]{article}
\usepackage[textwidth=6.5in,textheight=8.5in,top=1.25in]{geometry}
\newtheorem{definition}{Definition}

\usepackage{endnotes}
\let\footnote=\endnote

\usepackage{pst-all}
\usepackage{subfigure}
\usepackage{url}
\usepackage{graphics}
\usepackage{amsmath}
\usepackage{amssymb}
\usepackage{rotating}
\usepackage{booktabs}
\usepackage{multirow}
\usepackage{algorithm}
\usepackage{algorithmic}
\usepackage{tikz}

\DeclareMathOperator*{\argmax}{arg\,max}

\usepackage{natbib}
\usepackage[utf8]{inputenc}

 \bibpunct[, ]{(}{)}{,}{a}{}{,}%
\begin{document}

\title{A Primal Decomposition Algorithm for the Two-dimensional Bin Packing Problem}
\author{Jean-Fran\c{c}ois C\^ot{\'e}$^{(1)}$, Mohamed Haouari$^{(2)}$, Manuel Iori$^{(3)}$\\ \\
 $(1)$ CIRRELT, Universit\'e Laval,\\
 2325, rue de la Terrasse, Qu\'ebec, Qu\'ebec, Canada G1V 0A6\\
{\tt jean-francois.cote@fsa.ulaval.ca}\\
$(2)$ DMIE, College of Engineering, Qatar University, Doha, Qatar\\ 
{\tt    mohamed.haouari@qu.edu.qa}\\
$(3)$ DISMI, University of Modena and Reggio Emilia,\\ Via Amendola 2, 42122 Reggio Emilia, Italy\\
{\tt manuel.iori@unimore.it}
}
\maketitle

\vspace{-0.15 in}

\begin{abstract}{
The Two-dimensional Bin Packing Problem calls for packing a set of rectangular items into a minimal set of larger rectangular bins. Items must be packed with their edges parallel to the borders of the bins, cannot be rotated and cannot overlap among them. The problem is of interest because it models many real-world applications, including production, warehouse management and transportation. It is, unfortunately, very difficult, and instances with just 40 items are unsolved to proven optimality, despite many attempts, since the 1990s. In this paper, we solve the problem with a  combinatorial Benders decomposition that is based on a simple model in which the two-dimensional items and bins are just represented by their areas, and infeasible packings are imposed by means of exponentially-many no-good cuts. The basic decomposition scheme is quite naive, but we enrich it with a number of preprocessing techniques, valid inequalities, lower bounding methods, and enhanced algorithms to produce the strongest possible cuts. The resulting algorithm behaved very well on the benchmark sets of instances,  improving on average upon previous  algorithms from the literature and solving for the first time a number of open instances.

\vspace{0.1 in}
\noindent \textbf{Keywords}: two-dimensional bin packing problem; exact algorithm; {Benders} decomposition; combinatorial {Benders} cut
}
\end{abstract}

\section{Introduction}\label{sec:introduction}

In the \emph{Two-dimensional Bin Packing Problem} (2D-BPP),  we are given a set of rectangular items and a large number of rectangular bins. The aim is to pack the items {into} the minimum number of bins, in such a way that items are not rotated, do not {overlap}, and are packed with their edges parallel to the borders of the bins.

The 2D-BPP belongs to the area of \emph{Cutting \& Packing} (C\&P) problems, and has been intensively studied since the 1990s \citep{MV98} because it can model a large number of real-world applications. In packing applications, it serves, for instance, in determining the minimum number of cells to accommodate items in a warehouse, or in finding the minimum number of containers to load items to be shipped. In cutting applications, it is useful {for computing} the minimum number of plates required to produce a set of demanded items in productions involving materials such as steel, wood and glass. In addition, it has applications in scheduling theory, when assigning jobs with given length and resource consumption, in telecommunications, when allocating tasks involving processing power and delays, and in vehicle routing, when determining the minimum number of vehicles that can load a set of items that cannot be stacked one over the other. We refer the reader interested in the C\&P literature to the typology by \citet{WHS07}, to the comprehensive book by \citet{S18}, and to the {survey} by \citet{LMMV14}. In addition, for the standard \emph{one-dimensional Bin Packing Problem} (1D-BPP) we refer to the surveys by \citet{V02} and \citet{DIM16,DIM18}. For interesting problem generalizations we refer, among others, to the works by \citet{BW13} and \citet{CPT12} (multi-dimensional packing), by \citet{IM10,IM13} (integrated routing and packing problems), by \citet{CO11} and \citet{NDIS18} (scheduling and packing problems), by \citet{MAJ18} (lot-sizing and cutting stock problems), and by \citet{TP18} (load balancing and packing).

Despite the large number of heuristic and exact algorithms developed for its solution, the 2D-BPP remains a very challenging problem in practice, as instances with just 40 items  are still unsolved to proven optimality. The state-of-the-art exact method \citep{PS07} is based on a dual decomposition of the 2D-BPP. In a master problem, the 2D-BPP is modeled as a set covering problem, where each column represents a feasible 2D packing of a bin comprising one or more items. {Since the number of such packings is} exponentially large, the authors implemented a column generation algorithm that first solves the master with a reduced set of columns, then determines if columns with negative reduced cost exists and, if any is found, adds  them to the master. The process is reiterated until a set of packings comprising all items and requiring a minimum number of bins is found.

In this paper, we follow a different approach. We invoke a primal decomposition method in which the 2D-BPP is modeled in a master problem by means of a descriptive model in which each variable states whether an item is packed into a given bin or not. We do not take into consideration the actual size of the items and of the bins, but instead ensure that the sum of the item areas does not exceed the area of the bin in which they are packed. The master we use thus follows the footsteps of the model proposed by \citet{MT90} for the 1D-BPP, which is itself based on the seminal work by \citet{K60}. Whenever an integer solution is found, we check in a sub-problem if the packings of the bins are 2D-feasible, that is, we check the existence of feasible 2D packings without overlapping. If such packings are obtained for all bins, then the solution found is proven optimal. Otherwise, no-good cuts are added to the master problem, which is then solved again.

The basic decomposition just described is particularly slow in practice, because: (1) the master problem has a very weak continuous relaxation and is subject to large symmetries; (2) the sub-problem is not only NP-hard but also very difficult in practice; (3) the classical no-good cuts that are iteratively added to the master problem can only remove very limited portions of the search space. We consequently developed enhanced techniques to take care of each of these drawbacks: (1) valid inequalities and preprocessing techniques are adopted to improve the continuous relaxation of the master problem; (2) a tailored exact algorithm is developed for a fast solution of the sub-problem; (3) several ways to improve the no-good cuts are invoked, so as to strengthen them as much as possible. The resulting algorithm is very effective and can solve a number of open benchmark instances from the literature, comparing favorably with the previous solution approaches for the 2D-BPP.

The idea that we adopted derives from recent researches that obtained good computational results on a number of difficult C\&P problems.
\citet{CM04} and \citet{BB07} solved the \emph{Two-dimensional Knapsack Problem} (2D-KP), that is, the problem of packing a set of valued rectangular items into a single rectangular bin by maximizing the total value. In both approaches, the master is a mathematical model for the one-dimensional knapsack problem, whereas the the sub-problem consists of a 2D feasibility check. \citet{CDI14} used a primal decomposition method to solve the {\em Two-Dimensional Strip Packing Problem} (2D-SPP), the problem of orthogonally packing a given set of items without overlapping in a strip of given width and infinite height, by minimizing the height used for the packing. Their approach is different from the one presented in this paper as it involves solving in the master problem a relaxation in which the items can be cut into vertical slices that are then packed contiguously in the strip, and then checking in the sub-problem if all slices of an item can be packed at the same vertical height, for each item, thus building a feasible packing, if any. Their idea is indeed at the basis of the exact algorithm that we adopted to solve our sub-problem. A similar approach was later adopted by \cite{CGP14} for the 2D-SPP with unloading constraints, and by \citet{DIM17} for the 2D-SPP with item rotation. Very recently, \citet{DDIM19} developed a decomposition for the \emph{Multiple Knapsack Problem}, that makes iterated use of arc-flow models \citep{V99}, and in particular of the reflect model by \citep{DI19}. The first time the reflect model is used in a master problem to select the subset of items of maximum profit, and the second time it is used in a sub-problem to check if these items can be partitioned into the multiple knapsacks.

We believe this paper has a number of interesting contributions: (i) we develop a new exact algorithm for the 2D-BPP that adopts a different approach from the ones available in the literature; (ii) we gather together state-of-the-art techniques to effectively tackle C\&P problems, including combinatorial  cuts \citep{CF06}, preprocessing techniques \citep{BM10}, dual feasible functions \citep{ACCR16}, conservative scales \citep{BKRS13}, valid inequalities and lifting techniques \citep{KL10}; (iii) we solve for the first time a number of open instances from the literature; (iv) we show how the developed ideas can be adjusted to solve other difficult well-known C\&P problems.

The remainder of this paper is organized as follows. Section \ref{sec:2DBPP} formally defines the problem and presents the mathematical model that is at the basis of our decomposition. Section \ref{sec:solution-method} presents our overall solution method. Sections \ref{sec:preprocessing}, \ref{sec:lower_bounds} and \ref{sec:valid-inequalities} list, respectively, some preprocessing techniques, lower bounding algorithms and valid inequalities that are useful in speeding up the solution process. Section \ref{sec:combinatorial-cuts} describes the key sub-problem that we need to tackle iteratively in the decomposition and shows the main algorithms that we developed for producing enhanced cuts. Extensive computational results on benchmark instances from the literature and comparison with the state-of-the-art algorithms are provided in Section \ref{sec:results}. Despite the many efforts, instances with just 40 items remain unsolved to proven optimality, and the situation is even worse in other related 2D C\&P problems. Relevant hints for future research directions, both on the 2D-BPP and on other problems, are thus provided in the concluding Section \ref{sec:conclusions}.

\section{Problem Description and Mathematical Formulation}\label{sec:2DBPP}

In the 2D-BPP,  we are given a set $N = \{1, 2, \dots, n\}$ of rectangular items of width $w_j$, height $h_j$ and area $a_j = w_j h_j$, $j \in N$, and a set $B = \{1, 2, \dots, m\}$ of rectangular bins all having width $W$, height $H$ and area $A = WH$. The aim is to pack all items into the minimum number of bins without overlapping and by ensuring that the item edges are parallel to the borders of the bins. Rotating items by 90 degrees is not allowed. In the following, we suppose that $m$ is large enough to allow a feasible solution for the problem, i.e., $m$ takes the value of a valid upper bound on the optimal solution value. 

A number of interesting algorithms have been developed to solve the 2D-BPP. In terms of exact methods, we mention the seminal branch-and-bound  by \citet{MV98}, the dual decomposition method by {\citet{PS07}}, the method based on the iterative decomposition of the set of items into two disjoint subsets by \citet{CCM07} and the exact enumeration scheme for general dimensions C\&P problems by \citet{FSV07}. In terms of (meta)heuristics methods, we mention the Guided Local Search by \citet{FPZ03}, the Tabu Search by \citet{LMV04}, the set-covering-based heuristic by \citet{MT06}, the combination of Greedy Randomized Adaptive Search with Variable Neighborhood Descent by \citet{PAOT10} and {the goal-driven metaheuristic by \citet{WOZL13}}. In terms of lower bounds, a recent theoretical and experimental study of fast algorithms has been presented by \citet{SH18}.

Let us introduce two families of binary variables: $y_i$ takes the value 1 if bin $i$ is open, for $i \in B$;   $x_{ij}$ takes the value 1 if item $j$ is assigned to bin $i$, for $i \in B$, $j \in {N}$. Let us also define $S \subseteq N$ as a generic subset of items, and $\mathcal{S} \subseteq 2^N$ as the class of subsets of items that cannot be feasibly packed into a single bin. The 2D-BPP can then be modeled as the following \emph{Integer Linear Program} (ILP):
\begin{eqnarray}
 \mbox{(2D-BPP)} \qquad \min z = \sum_{i \in B} y_i, \label{form1:obj}\\
 \sum_{i \in B} x_{ij} = 1 &&  j \in N, \label{form1:eq1}\\
 \sum_{j \in N} a_j x_{ij} \leq A y_i && i \in B , \label{form1:eq2}\\
 \sum_{j \in S} x_{ij} \leq |S| - 1  && i \in B , S \in \mathcal{S} , \label{form1:eq3}\\
  x_{ij} \in \{0,1\} && i \in B, j \in {N}, \label{form1:eq4}\\
  y_i \in \{0,1\} && i \in B. \label{form1:eq5}
\end{eqnarray}
The objective function \eqref{form1:obj} asks to minimize the number of open bins. Constraints \eqref{form1:eq1} impose that each item is packed {into} a bin. Constraints \eqref{form1:eq2} ensures that only open bins are used and that the sum of the item areas assigned to a bin does not exceed the bin area. Constraints \eqref{form1:eq3} are the classical no-good cuts that forbid assigning more than $|S|-1$ items to a bin for those subsets $S \in \mathcal{S}$. Constraints \eqref{form1:eq4} and \eqref{form1:eq5} define the variable domains.
Concerning constrains \eqref{form1:eq3}, it is worth noting {that: (i) they are exponentially many and it is better to include only cuts corresponding to minimal infeasible subsets; and (ii)} computing if a given subset $S$ belongs to class $\mathcal{S}$ is NP-complete, as it corresponds to the classical \emph{Two-dimensional Orthogonal Packing Problem} (2D-OPP), see, e.g., \citet{CCM07b, CJCM08}.
 
Model \eqref{form1:obj}--\eqref{form1:eq2}, \eqref{form1:eq4}, \eqref{form1:eq5} is the classical ILP of the 1D-BPP by \citet{MT90}. It is known to be weak because of the large number of symmetries (a solution can be transformed into an equivalent one just by swapping items in a bin with items in another bin), and of the weak continuous relaxation. Indeed, let $c(\mbox{1D-BPP})$ be the optimal value of the 1D-BPP model relaxation obtained by replacing \eqref{form1:eq4} with $0 \leq x_{ij} \leq 1$ and \eqref{form1:eq5} with $0 \leq y_j \leq 1$, and let $z(\mbox{1D-BPP})$ be the optimal solution value. It is known that $c(\mbox{1D-BPP})/z(\mbox{1D-BPP})$ tends to 1/2 and the ratio is asymptotically tight (consider an instance with $n$ items of area $A/2+\epsilon$: $z(\mbox{1D-BPP})=n$, while $c(\mbox{1D-BPP})=n/2+1$).

Also the continuous relaxation of the entire 2D-BPP model can be quite far from its optimal solution. Let $c(\mbox{2D-BPP})$ and $z(\mbox{2D-BPP})$ be the continuous relaxation value and the optimal value, respectively, of  model \eqref{form1:obj}--\eqref{form1:eq5}. We can notice that $c(\mbox{2D-BPP})/z(\mbox{2D-BPP})$ can be as {bad} as 1/3, by considering a simple instance of 3 items each having width $W/2+\epsilon$ and height $H/2+\epsilon$: we get $c(\mbox{2D-BPP})=1$ as three such items can be packed into the same bin while respecting \eqref{form1:eq2}, while $z(\mbox{2D-BPP})=3$. \citet{MV98} used a generalization of this simple instance to prove that the ratio can be as bad as 1/4 and that this is asymptotically tight.

Despite the symmetries and the weak relaxation, it is reasonable to employ model \eqref{form1:obj}--\eqref{form1:eq2}, \eqref{form1:eq4}, \eqref{form1:eq5} as a starting tool to solve the 2D-BPP. In the experiments in \citet{DIM16, DIM18}, this model, executed for a minute on a standard PC using Cplex 12.6.0 as ILP solver, could solve to optimality the majority of {1D-BPP} instances with up to 200 items, and a number of larger ones with up to 1000 items. This is fair enough when addressing the 2D-BPP, whose benchmark instances, some of which are still unsolved to proven optimality, contain from 20 to 100 items. 

\section{Overall solution algorithm}\label{sec:solution-method}

The pseudocode of the method that we developed to solve the 2D-BPP is shown in Algorithm \ref{alg:overall}. It basically consists in a \emph{Branch-and-Cut} (B\&Cut) that solves model \eqref{form1:obj}--\eqref{form1:eq5} by adding constraints of type \eqref{form1:eq3} only when needed. The pseudocode already contains a description of the main steps of the algorithm. We mention a number of additional relevant details: 
\begin{itemize}
\item To compute the initial solution we use the heuristic by \citet{PAOT10}. In the computational tests, we will asses the impact of this algorithmic component by also attempting a test in which $U_0$ is set  to the best known value in the literature for the instance; 
\item Differently from other B\&Cut algorithm, we only separate inequalities at integer nodes. This might slow down the increase in the overall lower bound, but has the clear advantage of minimizing the number of 2D-OPP checks, a very difficult and time consuming component of the algorithm; 
\item Any time a 2D-OPP check is performed, the corresponding set $S$ is stored in a long term memory using a simple hash function, so as to avoid redundant 2D-OPP checks;
\item Still with the aim of minimizing the number of 2D-OPP checks performed, {if a proven-infeasible packing is found for a bin, then the successive 2D-OPP checks are not performed if they involve $\tilde n$ or more items, with $\tilde n = 18$ in our experiments}.
\end{itemize}

\begin{algorithm} [htb]
\caption{Solution algorithm for the 2D-BPP} \label{alg:overall}
\begin{algorithmic}
\REQUIRE $N$: set of items, $B$: set of bins
\STATE Reduce the instance with the preprocessing techniques of Section \ref{sec:preprocessing}
\STATE Compute an initial solution of value $U_0$ with a heuristic from the literature
\STATE Compute an initial lower bound $L_0$ with the methods from Section \ref{sec:lower_bounds}
\STATE {\bf if} $L_0=U_0$ {\bf then} {\bf return} the optimal solution
\STATE Initialize model \eqref{form1:obj}--\eqref{form1:eq5} with no cut of type \eqref{form1:eq3} but all inequalities from Section \ref{sec:valid-inequalities}
\STATE Solve model model \eqref{form1:obj}--\eqref{form1:eq5} in a B\&Cut fashion
\STATE \hspace*{0.5cm} \textbf{for} {each integer solution found in the B\&Cut search}
\STATE \hspace*{1cm} Solve a 2D-OPP for each packing of a bin with the method from Section \ref{subsec:OPP}
\STATE \hspace*{1cm} {\bf if}  all 2D-OPP checks return feasible, {\bf then} update the incumbent
\STATE \hspace*{1cm} {{\bf else} compute the corresponding no-good cuts of type \eqref{form1:eq3}, lift them with the techniques 
\STATE \hspace*{1.5cm} from Sections \ref{subsec:MIS} and \ref{subsec:lifting} and then add them to the model}
\STATE \hspace*{0.5cm} \textbf{end-for}
\STATE \textbf{return} the incumbent solution
\end{algorithmic}
\end{algorithm}

\section{Preprocessing}\label{sec:preprocessing}

In this section, we describe the techniques that we developed to shrink the sizes of the bins and enlarge the sizes of the items (Section \ref{subsec:preproc-shrink}) and preliminary packing some items (Section \ref{subsec:preproc-pack}). All such techniques have the aim of reducing the complexity of an instance so as to obtain problems that are easier to be solved in practice.

\subsection{Shrinking the bins and enlarging the items}\label{subsec:preproc-shrink}
As shown by \citet{APT09} for the 2D-SPP, if  there are no combinations of item widths whose sum is equal to $W$, then it is possible to conveniently reduce $W$ without losing any optimal solution. To this aim, we compute by dynamic programming the maximum value $W^* \leq W$ that can be produced by combining the item widths as
\begin{eqnarray}
\label{eq:min_width}	W^* &=& \textstyle \max \left \{z=\sum_{j \in N}z_j w_j : z \leq W,\ z_j \in  \{0,1\}, j \in N \right \} 
\end{eqnarray}
and then we set $W=W^*$. For the 2D-BPP, the same reasoning can be applied to $H$ as well, so we use a variant of \eqref{eq:min_width} in which $w_j$ and $W$ are replaced by $h_j$ and $H$, respectively, to compute the maximum usable height $H^*$, and then set $H=H^*$.

A similar procedure can then be applied to increase the item sizes. This was first noticed for the 2D-KP by \citet{BHM02} and then applied to the 2D-BPP by \citet{CCM07}. For any item $j \in N$, we compute by dynamic programming the values
\begin{eqnarray}
\label{eq:min_item_width}	w^*_j &=& \textstyle W - \max \left \{z=\sum_{k \in N\backslash \{j\}}z_k w_k : z \leq W-w_j,\ z_k \in  \{0,1\}, k \in N\backslash \{j\}\right \} 
\end{eqnarray}
and then use them to lift the original sizes by setting $w_j=w^*_j$. Similarly, a variant of \eqref{eq:min_item_width} involving heights is used to compute the maximal heights $h^*_j$ for any $j \in N$.

We propose a simple improvement on the above preprocessing techniques that takes into account the fact that, as discussed in Section \ref{sec:valid-inequalities} below, some items cam be forbidden to be packed in some bins without losing optimality.
Let $N_i$ be the set of items that are allowed to be packed inside a bin $i \in B$. We first shrink the width of $i$ by computing
\begin{eqnarray}
\label{eq:min_width2}	W^*_i &=& \textstyle \max \left \{z=\sum_{j \in N_i} z_j w_j : z \leq W,\ z_j \in  \{0,1\}, j \in N_i \right \} 
\end{eqnarray}
and then lift the widths of the items in $N_i$ to 
\begin{eqnarray}
\label{eq:min_item_width2} \hspace*{-0.3cm} w^*_{ij} &=& \textstyle W - \max \left \{z=\sum_{k \in N_i\setminus \{j\}} z_k w_{ik}: z \leq W_i-w_{ij},\ z_k \in \{0,1\}, k \in N_i\setminus \{j\} \right \}. 
\end{eqnarray}
A similar procedure involving heights is adopted to produce the minimal height $H^*_i$ of a bin and the maximal heights $h^*_{ij}$ of the items.

{For simplicity, during the computation of our lower bounds (Section \ref{sec:lower_bounds}) we made use only of the preprocessing \eqref{eq:min_width} and \eqref{eq:min_item_width}.  
We used instead \eqref{eq:min_width2} and \eqref{eq:min_item_width2}, and the corresponding variants involving heights, in the mathematical model and in the procedure developed for its solution via B\&Cut.}

\subsection{Packing and removing some items}\label{subsec:preproc-pack}

Starting from \citet{MMV03}, a number of authors proposed ideas aimed at preliminary packing some items of a 2D-SPP instance at the bottom of the strip without losing optimality. These ideas were later extended to the 2D-BPP by \citet{CCM07}. Here, we replicate the preprocessing techniques in Sections 3.1 and 3.2 of \citet{CCM07}. We describe our implementation in detail so as to allow its replicability.

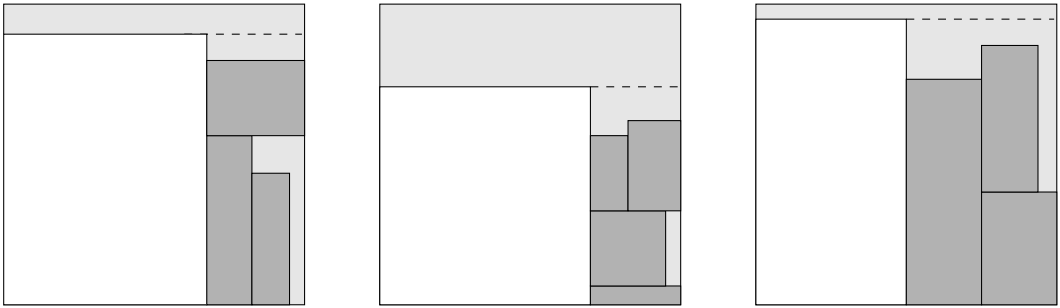
\begin{figure}[ht]
\begin{centering}
\begin{tikzpicture}[shorten >=1pt,scale=1.0,every node/.style={scale=1.0}]
\tikzstyle{vertex}=[circle,fill=black!25,minimum size=14pt,inner sep=0pt]
   \filldraw[fill=gray!20!white, draw=black] (0,0) rectangle (4,4);
   \filldraw[fill=gray!20!white, draw=black] (5,0) rectangle (9,4);
   \filldraw[fill=gray!20!white, draw=black] (10,0) rectangle (14,4);
   \draw[dashed] (2.4,3.6) -- (4,3.6) ;
   \draw[dashed] (7.8,2.9) -- (9,2.9) ;
   \draw[dashed] (12,3.8) -- (14,3.8) ;

   \foreach \x/\y/\w/\h in {0/0/2.7/3.6, 5/0/2.8/2.9, 10/0/2/3.8}
      \filldraw[fill=white, draw=black] (\x,\y) rectangle (\x+\w,\y+\h);

   \foreach \x/\y/\w/\h in {2.7/0/0.6/2.25, 2.7/2.25/1.3/1, 3.3/0/0.5/1.75,  7.8/0/1.2/0.25, 7.8/0.25/1/1, 7.8/1.25/0.5/1,8.3/1.25/0.7/1.2, 12/0/1/3, 13/0/1/1.5, 13/1.5/0.75/1.95 }
      \filldraw[fill=gray!60!white, draw=black] (\x,\y) rectangle (\x+\w,\y+\h);

\end{tikzpicture}
\end{centering}
\caption{Fixing some items by preprocessing}
\label{fig:fixing_items}
\end{figure}

Let $C$ be a subset of items whose width is greater than $W/2$. Let $w^* = W-\min\{w_i: i \in C\}$ and $R = \{j \in N : w_j \leq w^*\}$. Items in $R$ can only be packed side by side with the items in $C$ or in some other bins. If they all fit side by side with the items in $C$, then the items in $R$ can be all removed from the instance and the items in $C$ can be all enlarged to have width $W$. Figure \ref{fig:fixing_items} shows an example. The white items are those in $C$ and the gray items are those in $R$. In the example, items in $R$ can be removed because they all fit below the dashed lines.

To test if such a feasible packing exists for a given set $C$, we first build a set $B'$ of $|C|$ bins. Each bin has width $W-w_i$ and height $h_i$, for $i \in C$. Then, we create the corresponding set $R$ and  invoke a heuristic to look for a feasible packing of the items in $R$ into the bins of $B'$. In our implementation, we use a heuristic that works as follows. We sort bins from the one with smallest width to the one with largest width. We fill one bin at a time until no more items can enter it or all items have been packed. To fill a bin, we invoke the heuristic by  \citet{LZS11}. If this heuristic finds a feasible packing, then the items in $R$ are removed from the instance and the widths of all items in $C$ are set to $W$, otherwise we proceed with a new tentative set $C$.

To scan for the relevant sets $C$, in our implementation we first order the items by non increasing width, breaking ties by non increasing height, and then perform two loops. In the first loop, we check all sets $C$ having size one. In the second loop, $C$ is initially composed by the first two largest items. If a feasible packing is found, then the items in $C$ and in the corresponding set $R$ are temporarily removed from the instance and the algorithm restarts by initializing $C$ to the next two largest items. Otherwise, the next largest item is added to the current $C$ and the algorithm continues checking for feasible packings until the next largest item has a width smaller than or equal to $W/2$.

Following \citet{CCM07}, the above technique can be adapted to deal with two other cases. Firstly, we consider tall items having height $h_i > H/2$. Let $h^* = H-\min\{h_i : i \in C\}$ and $R = \{j \in N : h_j \leq h^*\}$. If all items in $R$ fit inside the bins $B' = \{ (w_i, H-h_i) : i \in C\}$, then they are removed from the instance and the height of each item in $C$ is set to $H$. The implementation is the same as the one used for the large items, but considers widths instead of heights when scanning the tentative sets $C$.
Secondly, we consider all items having both $w_i > W/2$ and $h_i > H/2$. We create a set $C$ and let $R = \{j \in N : \exists i \in C \mbox{ such that } w_j + w_i \leq W \mbox{ or } h_j + h_i \leq H \}$. If a feasible packing of the items in $R \cup C$ inside $|C|$ bins of dimensions $(W,H)$ is found, then all items of $R$ are removed from the instance and the dimensions of the items of $C$ are enlarged to $W$ and $H$. The implementation is the same as the one used for the large items, but considers areas instead of widths when scanning the tentative sets $C$.

The above preprocessing methods are run in loop until no improvement is found. At the end, we are typically left with an instance that contains a reduced number of small items and a number of large items with even larger dimensions. All items having dimensions $(W,H)$ are packed {into} separate bins and then removed from the instance.

\section{Lower Bounds}\label{sec:lower_bounds}

During the last decades a large number of lower bounding techniques have been produced for the 2D-BPP (and for C\&P problems in general). The most recent and effective ones are based on the concept of dual feasible functions and discrete dual feasible functions.
\begin{definition} A function $f:[0,1] \rightarrow [0,1]$ is a \emph{dual feasible function} (DFF), if for any finite set $S$ of  positive real numbers  the following relation holds
\begin{equation}
	\sum_{x \in S} x \leq 1 \Rightarrow \sum_{x \in S} f(x) \leq 1.
\end{equation}
\end{definition}
\begin{definition} A function $f:[0,C] \rightarrow [0,C']$ is a discrete DFF, if for any finite set $S$ of integer numbers the following relation holds
\begin{equation}
	\sum_{x \in S} x \leq C \Rightarrow \sum_{x \in S} f(x) \leq f(C) = C'.
\end{equation}
\end{definition}
\citet{CAV10} highlight that for any DFF there is a corresponding discrete DFF and {vice-versa}, so it is equivalent to use one definition or the other. In our work, we consider the discrete DFF. Roughly speaking, discrete DFF are used to modify the widths and/or the heights of items and bins. New areas are computed by attempting one or more combinations of modified widths and heights, and then used to compute new continuous area bounds.

Our lower bound $L_0$ is based on two components, both making use of discrete DFF. The first one is the lower bound $L^2_{CCM}$ by \citet{CCM07}, that is based on the use of three different discrete DFF. Let $k$ be an integer parameter taking any value in the interval $[1, C/2]$. The first DFF is simply: 
\begin{alignat}{1}
	f_0^k(x) =  \begin{cases}
	C & \mbox{ if } x > C-k,\\
	x & \mbox{ if } C-k \geq x \geq k,\\
	0 & \mbox{ otherwise}.
	\end{cases}
\end{alignat}
The second is a discrete \emph{Data-dependent DFF}, that is, a discrete DFF that works for a particular instance and not necessarily for any general combination of numbers (we refer to \citealt{CAV10}  for further details). Let $J = \{j \in N: C/2 \leq x_j \leq k\}$ (where $x_j$ and $C$ can be either widths or heights of items and bin). Let $CKP(C,J)$ be a \emph{Cardinality Knapsack Problem} that aims at maximizing the number of items from $J$ that can be packed {into} a bin of size $C$, and let $z(CKP(C,J))$ be its optimal solution value. The second function is stated as
\begin{alignat}{1}
	f_1^k(x) =  \begin{cases}
	z(CKP(C,J)) - z(CKP(C-x,J)) & \mbox{ if } x > C/2,\\
	1 & \mbox{ if } C/2 \geq x \geq k,\\
	0 & \mbox{ otherwise}.
	\end{cases}
\end{alignat}
The third function is obtained by setting
\begin{alignat}{1}
	f_2^k(x) =  \begin{cases}
	2(\lfloor \frac{C}{k} \rfloor - \lfloor \frac{C-x}{k} \rfloor), & \mbox{ if } x > \frac{C}{2},\\
	\lfloor \frac{C}{k} \rfloor, & \mbox{ if } x = \frac{C}{2},\\
	2\lfloor \frac{x}{k} \rfloor, & \mbox{ if } x < \frac{C}{2}.
	\end{cases}
\end{alignat}
The three functions are computed for all possible parameter values and combined together in all possible ways, to obtain the highest bound of type
\begin{equation}\label{eq:CMM}
L^2_{CCM} =	\textstyle \max_{u, v, k, l} \left \lceil \sum_{j \in N} (f_u^k(w_j) f_v^l(h_j)) / (f_u^k(W) f_v^l(H)) \right \rceil,
\end{equation}
where $u, v = 0, 1, 2$ and $k, \l \in [1, C/2]$. For the details on the computation of $L^2_{CCM}$ we refer to Section 4 of \citet{CCM07}.

The lower bound in \eqref{eq:CMM} is fast as can be computed in polynomial time, and, as shown by \citet{SH18}, is the one providing the best value among all the 2D-BPP fast lower bounds on all benchmark instances. To look for improved values, in case $L^2_{CCM}$ is not enough to prove optimality of a solution, our second $L_0$ component is thus based on a more complex {non-polynomial} procedure.

The problem of finding the best possible lower bound derived by discrete DFF can be formulated as the following bilinear program with disjoint constraints:
\begin{eqnarray}
 \max  \sum_{j \in N} \tilde{w_j} \tilde{h_j}, \label{lb:dff_obj} \\
 \sum_{j \in S} \tilde{w_j} \leq W &&  S \in \mathcal{S}^W, \label{lb:dff_const1}\\
 \sum_{j \in S} \tilde{h_j} \leq H &&  S \in \mathcal{S}^H, \label{lb:dff_const2}\\
\tilde{w_j}, \tilde{h_j} \geq 0 && j \in N \label{lb:domain},
\end{eqnarray}
where the decision variables $\tilde{w_j}$ and $\tilde{h_j}$ represent the new sizes of item $j$, $\mathcal{S}^W = \{S \subseteq N: \sum_{j \in S} w_j \leq W\}$ and  $\mathcal{S}^H = \{S \subseteq N : \sum_{j \in S} h_j \leq H\}$. The objective function \eqref{lb:dff_obj} maximizes the sum of the modified item areas, in such way that the items in any feasible subset can still be be packed side by side, as stated by \eqref{lb:dff_const1}, or one above the other, as stated by \eqref{lb:dff_const2}.

\citet{CM09} investigated ways to solve exactly this model to obtain the best possible bounds. They were able to prove optimality for several open instances, however, their computing times were very high on many instances. To tackle this difficulty, \citet{BKRS13} (who use the name of \emph{conservative scales} to denote $\tilde{w_j}$ and $\tilde{h_j}$) proposed a heuristic that iteratively fixes the variables of one of the two dimensions and solves the reduced problem on the remaining dimension. In this way, the difficult bilinear program reduces to a series of simple linear programs. 

In our work, we implemented the same procedure of \citet{BKRS13} and obtained a lower bound that we call $L_{BKRS}$. 
We perform a series of $\eta$ iterations. At each iteration $k = 1, \dots, \eta$, we obtain new vectors of modified item widths $w^k$ and item heights $h^k$, by solving the two independent linear programs
\begin{eqnarray}
 w^k = \textstyle \argmax_{\tilde{w}} \left \{ \sum_{j \in N} \tilde{w_j} h_j^{k-1} : \sum_{j \in S} \tilde{w_j} \leq W, S \in \mathcal{S}^W \right \}, \label{cs:prob1} \\
 h^k = \textstyle \argmax_{\tilde{h}} \left \{ \sum_{j \in N} \tilde{h_j} w_j^{k-1} : \sum_{j \in S} \tilde{h_j} \leq H, S \in \mathcal{S}^H \right \}. \label{cs:prob2}
\end{eqnarray}
The values of $w_j^{k-1}$ and $h_j^{k-1}$ are obtained by the previous iteration, and initialized to, respectively, $w^0=w$ and $h^0=h$ (i.e., the original item sizes). 
An anti-stalling procedure is invoked in case the modified widths and heights do not change from one iteration to the next. This is described in detail in Section 4.3 of \citet{BKRS13}.
Once all iterations have been performed, we obtain $\eta^2$ lower bounds by attempting any combination of $w_k$ and $h_l$ for $k, l = 1, \dots, \eta$, namely
\begin{equation}\label{eq:BKRS}
L_{BKRS} =	\textstyle \max_{k, l \in 0, \dots, \eta} \left \lceil \sum_{j \in N} (w_j^k h_j^l) / (W^k H^l) \right \rceil .
\end{equation}
We then set $L_0 = \max \{L^2_{CCM}; L_{BKRS}\}$.

\section{Valid Inequalites}\label{sec:valid-inequalities}

To speed up the convergence of the 2D-BPP model, we implemented some families of simple and enhanced valid inequalities that are added at the root node of the model. 

\subsection{Simple Inequalites}\label{subsec:simple-inequalities}

A classical way to reduce symmetry in the 2D-BPP model is to allow using bin $i+1$ only if bin $i$ is also used. This can be imposed by:  
\begin{eqnarray}
y_{i} \leq y_{i+1} && i=1, 2,\dots, m-1. \label{ineq:1}
\end{eqnarray}

A second classical way to reduce symmetry is to allow item $j$ to be in bin $i$ only if $i \leq j$. This is obtained by setting: 
\begin{eqnarray}
x_{ij} = 0 && i \in B, j \in {N},  i > j. \label{ineq:2}
\end{eqnarray}

In addition, let us define two items $j$ and $k$ to be \emph{incompatible} if $(w_j + w_k > W)$ and $(h_j + h_k > H)$ both hold.
Let $C \subseteq N$ be a subset of items being all pairwise incompatible. No pair of items belonging to $C$ can be feasibly assigned to the same bin, {so we make use of the following clique inequalities:
\begin{eqnarray}
\sum_{j \in C} x_{ij} \leq 1 && i \in B. \label{ineq:3}
\end{eqnarray}
In our implementation, we simply computed all cliques of size 2. The corresponding inequalities \eqref{ineq:3} were added to the model after being lifted with the procedure described in Section \ref{subsec:lifting} below.}
 
{Now, let $T \subseteq$ be the clique of maximal cardinality, and let} us reorder $N$ in such a way that the first $|T|$ items are those of $T$. According to \eqref{ineq:2}, item 1 is assigned to bin 1, item 2 to either bin 1 or 2, item 3 to either bin 1, 2, or 3, and so on. Then, items in $T \backslash \{1\}$ can be forbidden from bin 1, items in $T \backslash \{2\}$ from bin 2, items in $T \backslash \{3\}$ from bin 3 and so on for the first $|T|$ bins. At the end, each item in $T$ is assigned to a single bin, and this can be simply imposed as:
\begin{eqnarray}
x_{ii} = 1 && i \in {T}. \label{ineq:3}
\end{eqnarray}
To obtain the strongest reduction from \eqref{ineq:3}, we first build a conflict graph $G$ where each node represents an item and each edge represents an incompatibility between a pair of items, and then use it to compute $T$. The problem of determining a maximal clique $T$ in $G$ is NP-hard, but can be quickly solved to optimality for small graphs. In our implementation, we used the branch-and-bound by \citet{KJ07}. {We also set $x_{ij} = 0$ for each $j \in N, i \in T$ such that $j$ is incompatible with $i$.}

\subsection{Inequalites Based on Dual Feasible Functions}\label{subsec:DFF-inequalities}

Let ${N}_i \subseteq {N}$ be the set of items that are allowed to enter bin $i$, for $i \in B$, after the applications of \eqref{ineq:2} and \eqref{ineq:3}.
As described by \citet{CAV10}, any discrete DFF can be used to produce a valid inequality for an ILP. In our case, given two discrete DFF, $f^1$ and $f^2$, we create a valid inequality for model \eqref{form1:obj}--\eqref{form1:eq5} having the form
\begin{eqnarray}\label{eq:DFFs}
	\sum_{j \in N_i} f^1(w_{ij})f^2(h_{ij})x_{ij} \leq f^1(W_i)f^2(H_i)y_i &&  i \in B. \label{eq:dff_area1}
\end{eqnarray}

Note that the inequality takes into account the specific information computed by the preprocessing techniques \eqref{eq:min_width2} and \eqref{eq:min_item_width2} for any bin $i \in B$ (i.e., the subest $N_i$ and the dimensions of bin and items). The procedures implemented to obtain $L_0$ in Section \ref{sec:lower_bounds} can now be reused to produce valid inequalities of type \eqref{eq:dff_area1}, by taking into consideration the new specific information on each bin. We found convenient to include in the model the $\alpha$ pairs of different discrete DFF that lead to the highest total area in \eqref{eq:CMM} , and the $\beta$ pairs of different discrete DFF that  lead to the highest total area in \eqref{eq:BKRS}.

\section{Infeasible Subsets of Items and Combinatorial Cuts}\label{sec:combinatorial-cuts}

Each time the B\&Cut finds a new potential incumbent solution, feasibility of the packing of each bin needs to be checked. This is obtained by the algorithm in Section \ref{subsec:OPP}. If infeasible, a cut of type \eqref{form1:eq3} is added to the model (a cut for each bin $i \in B$). This cut can be weak, so we improve it in two ways: we first look for the minimal subset of items that causes infeasibility (Section \ref{subsec:MIS}) and then lift the cut by adding variables to its left-hand side (Section \ref{subsec:lifting}).

\subsection{Solving the 2D Orthogonal Packing Problem}\label{subsec:OPP}

The 2D-OPP is one of the most difficult problems in the C\&P field. Many attempts have been devoted to its solution, including Branch-and-Bound algorithms \citep{MV98, CCM07}, 
graph-theoretical models \citep{FSV07} and constraint {programming} \citep{PS07, CJCM08}. In our work, we opted to update the algorithm originally developed by \citet{CDI14}, which is based on a combinatorial {Benders} decomposition. In this sense, our overall algorithm for the 2D-BPP can be seen as a combinatorial {Benders} decomposition that invokes an inner combinatorial {Benders} decomposition for each 2D-OPP solution.

The algorithm in \citet{CDI14} has been developed for the 2D-SPP, so aims at minimizing the height used for a packing in a strip of given width. It performs a series of 2D-OPP checks between an upper and a lower bound. For each check, it divides each item into vertical slices and invokes an ILP to determine the $x$-coordinates in which the slices are packed, by ensuring that: (i) slices belonging to the same item are packed contiguously and (ii) the first slice of each item is packed in a normal pattern \citep{H72, CW77}. If a solution is found for the $x$-coordinates, the existence of feasible $y$-coordinates for all items is determined by means of a dedicated Branch-and-Bound algorithm. If it turns out that no such $y$-coordinates exist, a feasibility cut is produced, lifted in many ways and added to the original ILP, and then the process is iterated. If instead feasible $y$-coordinates are found, the process terminates with a proven-optimal solution. 

We performed two simple modifications of the original algorithm. The first one consists in performing only a single 2D-OPP check at height $H$. The second one consists in reducing the number of tentative $x$-coordinates for the items in the ILP by using the \emph{Meet-in-the-Middle} patterns \citep{CI18}. This is a small set of tentative positions for the items that still guarantees that an optimal solution exists but also reduces consistently the size of the standard set of normal patterns. 

\subsection{Finding Minimal Infeasible Subsets of Items}\label{subsec:MIS}

If the 2D-OPP check returned infeasible for a given set $S$, then a cut of type \eqref{form1:eq3} needs to be added to the model, one for each bin $i \in B$. This is commonly known in the literature as \emph{{Benders} feasibility cut}, or simply \emph{no-good cut}. It is known to be weak in practice when the set $S$ is large: a new linear solution which would be feasible for the cut can be produced by setting $x_{ij} = (|S|-1)/|S|$ for all $j \in S$, and this value can be close to 1, so close to the original integer point we intend to cut, when $|S|$ is large. To obtain new effective cuts, it is important to find the minimal source of infeasibility, that means for us the \emph{Minimal Infeasible Subset} (MIS) of items that still makes $S$ not packable {into} a single bin. Let $C$ be such a MIS, the resulting cut is 
\begin{eqnarray}
 \sum_{j \in C} x_{ij} \leq |C| - 1  && i \in B , C \in \mathcal{S} , \label{eq:comb-cut}
\end{eqnarray}
and is known in the literature as \emph{combinatorial {Benders} cut} \citep{CF06}.

Note that moving from a cut of type \eqref{form1:eq3} to a cut of type \eqref{eq:comb-cut} is a key component to obtain successful combinatorial decompositions, not only for C\&P problems (\citealt{CDI14, CGP14, DIM17}), but also for general optimization problems as in the \emph{Logic-based {Benders} Decomposition} by \citet{H07}. We also would like to note that this was already a key component in the dual decomposition method by \cite{PS07} for the 2D-BPP. 

Finding a MIS is NP-complete, so we content us with a heuristic approach. We first remove {an item a time from $S$}, from the one of smallest area to the one of largest area. For each new tentative set, we execute the 2D-OPP check with a limited time limit. If infeasibility is proven, we continue iterating by removing the next item. Otherwise, we stop and produce a cut for the last proven-infeasible subset. We follow the same approach starting from $S$ for other $\gamma$ times, but in these cases we follow randomly generated orders for removing the items. At the end, we produce a maximum of $\gamma+1$ cuts with different MIS. 

Each such cut has to be adapted to a given bin $i \in B$, by considering the subset $N_i$ of items that are allowed to enter such bin after the application of \eqref{ineq:2} and \eqref{ineq:3}. This is simply obtained in the iterative process that we just described by disregarding (i) the insertion of items that would be in contrast with \eqref{ineq:2} and (ii) the removal of items that would be in contrast with \eqref{ineq:3}.

\subsection{Lifting the Cut}\label{subsec:lifting}

Our last improvement algorithm attempts to \emph{lift} the inequality \eqref{eq:comb-cut} by adding variables (multiplied by a positive coefficient) to its left-hand side while leaving untouched its right-hand side. In other words, we look for an item $j^* \in N \setminus C$ that guarantees that at most $|C|-1$ items from the set $C \cup \{j^*\}$ can be feasibly packed {into} a bin. If such item is found, the process is sequentially repeated by looking for other additional items.

Formally, we look for a \emph{lifted cover inequality} \citep{B75b, W75, KL10} of the form
\begin{eqnarray}
 \sum_{j \in C} x_{ij} + \sum_{j \in N \backslash C} \alpha_j x_{ij} \leq |C|-1 && i \in B , C \in \mathcal{S}, \label{eq:lifted-cover}
\end{eqnarray}
where $\alpha_j$ are non-negative integer coefficients. 

To obtain the coefficients, we make use of Algorithm \ref{alg:liftedcover}, which we now explain in detail.
Let 2D-KP($S$, $j^*$) denote an instance of a 2D-KP whose aim is to pack a set $S$ of two-dimensional items, each having width $w_j$, height $h_j$ and profit $p_j$, $j \in S$, {into} a bin of width $W$ and height $H$, with the additional constraint that item $j^* \notin S$ is forced to be in the solution. Let $z(\mbox{2D-KP($S$, $j^*$)})$ denote the optimal solution value of the problem. Initially, consider $S = C$ and set $p_j = 1$ for all $j \in S$ and $p_{j^*}=0$. Then, $z(\mbox{2D-KP($S$, $j^*$)})$ gives the maximum number of items from $C$ that can be packed together with $j^*$ {into} a single bin. If this value is equal to $|C|-1$, then it is possible to pack $j^*$ with all but one items from $C$, so we set $\alpha_{j^*}=0$ in \eqref{eq:lifted-cover} otherwise we would be forced to increase the right-hand side of the cut. If instead it is lower than $|C|-1$, then its insertion in the bin would cause the removal of one or more items from $C$ in any feasible packing. In particular, it would cause the removal of at least $\alpha_{j^*} = |C| - 1 - z(\mbox{2D-KP($S$, $j^*$)})$ items. We can thus use this coefficient without increasing the right-hand side of the cut. To proceed with the sequential search, we then include $j^*$ in $S$ with profit $p_{i^*}=\alpha_{j^*}$ and iterate the process by looking for the next tentative additional item.

\begin{algorithm}[htb]
\caption{Calculating the lifting coefficients} \label{alg:liftedcover}
\begin{algorithmic}
\REQUIRE $N$: set of items, $C \subseteq N$: an infeasible subset of items
\STATE $S:= C$
\STATE $p_j:= 1$ for $j \in S$ and $p_j:= 0$ for $j \in N \setminus S$
\FOR{each $j^* \in N \setminus C$}
	\STATE $\alpha_{j^*}:= |C| - 1 - z(\mbox{2D-KP($S$, $j^*$)})$
	\IF  {$\alpha_{j^*} > 0$} \STATE {$S:= S \cup \{j^*\}$ and $p_{j^*}:= \alpha_{j^*}$} \ENDIF
\ENDFOR
\RETURN{$\alpha$}
\end{algorithmic}
\end{algorithm}

A full execution of Algorithm \ref{alg:liftedcover} would require to invoke $|N \setminus C|$ times an exact algorithm for the 2D-KP. As this can be very time consuming, we opted to invoke, instead, a relaxation of the 2D-KP. To this aim, we make use of the \emph{bar relaxation} by \citet{S99}. 

We define a \emph{pattern} $k$ as a binary array $a_{jk}$, with $a_{jk}=1$ if item $j$ is in the pattern and 0 otherwise.  A pattern $k$ is said to be $H$-feasible if $\sum_{j \in N} a_{jk} h_j \leq H$ and $W$-feasible if $\sum_{j \in N} a_{jk} w_j \leq W$. Let $K^H$ and $K^W$ denote, respectively, the sets of all $H$-feasible and $W$-feasible patterns. We consider a relaxation of the 2D-KP, called 2D-UKP, in which we do not require the selected rectangular items to be packed without overlapping in the bins, but only impose that each selected item $j$ is contained in at least $h_j$ $W$-feasible patterns and in at least $w_j$ $H$-feasible patterns.

Let the integer variable $x_k$, respectively $y_k$, denote the number of times a pattern $k \in K^H$, respectively $k \in K^W$, is selected. 
Let also $z_j$ be a binary variable indicating if item $j \in N$ is in the solution or not. Then, the 2D-UKP can be modeled as the following ILP:
\allowdisplaybreaks
\begin{eqnarray}
 \mbox{(2D-UKP)} \qquad \max \sum_{j \in {N}} p_j z_j, \label{eq:ukp-z}\\
 \sum_{k \in K^W} a_{jk} y_k \geq h_j z_j &&  j \in {N}, \label{eq:ukp-y}\\
 \sum_{k \in K^H} a_{jk} x_k \geq w_j z_j &&  j \in {N}, \label{eq:ukp-x}\\
 \sum_{k \in K^W} y_k \leq H, && \label{eq:ukp-H}\\
 \sum_{k \in K^H} x_k \leq W, && \label{eq:ukp-W}\\
  z_j  \in \{0,1\} && j \in  {N}, \label{eq:ukp-zbin}\\
  x_k \geq 0 \mbox{ and integer} && k \in K^H, \label{eq:ukp-xk}\\
  y_k \geq 0 \mbox{ and integer} && k \in K^W. \label{eq:ukp-yk}
\end{eqnarray}
Objective function \eqref{eq:ukp-z} asks for the maximization of the total profit. Constraints \eqref{eq:ukp-y} and \eqref{eq:ukp-x} ensure, respectively, that if item $j$ is selected, then it appears at least $h_j$ times in $W$-feasible patterns and $w_j$ times in $H$-feasible patterns. Constraints \eqref{eq:ukp-H} and \eqref{eq:ukp-W} impose, respectively, that at most $W$ $H$-feasible patterns and $H$ $W$-feasible patterns are used. 

Let $z(\mbox{2D-UKP})$ be the optimal solution value of the problem. Solving the 2D-UKP to optimality is challenging because the problem is NP-hard and, in addition, its formulation contains an exponential number of variables $x_k$ and $y_k$. For this reason, we consider the continuous relaxation of the 2D-UKP, obtained by removing integrality from constraints \eqref{eq:ukp-zbin}--\eqref{eq:ukp-yk}. We solve this relaxation with a standard column generation algorithm, in which we initialize the problem with a limited number of patterns, and then look for new useful patterns by solving one-dimensional 0-1 knapsack problems. Let $c(\mbox{2D-UKP})$ be the optimal value of the continuous relaxation. Being a maximization problem, we get 
$c(\mbox{2D-UKP}) \geq z(\mbox{2D-UKP}) \geq z(\mbox{2D-KP})$.

We use the above relaxation to update Algorithm \ref{alg:liftedcover} as follows. Let $c(\mbox{2D-UKP($S$, $j^*$)})$ be the optimal relaxation value of an instance composed by a set $S$ of items with profit 1 and an additional item $j^*$ with profit 0, in which $j^*$ is forced {into} the bin by adding the constraint $z_{j^*}=1$. The only change in the algorithm concerns the computation of the coefficients, that now we obtain by using $\alpha_{j^*}:= \max\{0; |C| - 1 - c(\mbox{2D-UKP($S$, $j^*$)}) \}$. The $\max$ {operator} in the equation follows from the fact that $c(\mbox{2D-UKP($S$, $j^*$)})$ can be greater than $|C|-1$.

As previously discussed for \eqref{eq:comb-cut}, also \eqref{eq:lifted-cover} can be made stronger by considering the subsets $N_i$ of items that can enter the different bins $i \in B$. In our implementation, this is obtained by imposing additional changes in the 2D-UKP formulation when computing the $c(\mbox{2D-UKP($S$, $j^*$)})$ values. In detail, (i) we simply disregard from $S$ items that cannot enter bin $i$ because of \eqref{ineq:2}, and (ii) we force items that should enter the bin because of \eqref{ineq:3} by setting each time a $z$ variable to one. 

\section{Computational Results}\label{sec:results}

The exact algorithm has been coded in C++ and tested on the 500 benchmark instances originally proposed by \citet{BW87} and \citet{MV98} and widely used in the 2D-BPP literature. The benchmark set is divided into 10 classes according to the way bin and items have been generated. It contains instances with 20, 40, 60, 80 and 100 items, all available at the web site maintained by the \cite{DEI19}. The web site also contains a list of known lower and upper bounds. We adopted Cplex 12.8.0 as ILP solver, keeping its default settings but imposing it to run on a single thread.  All computational tests have been executed on an Intel Gold 6148 Skylake 2.4 Ghz. Similarly to what done by \citet{PS07}, we let our exact algorithm run for one hour. {Concerning the number of valid inequalities of type \eqref{eq:DFFs} to be added to the model, we set $\alpha = \beta = 700$ on the basis of preliminary experiments performed with 900 seconds of time limit. Having a high number of such inequalities allows the MILP model to get a good initial lower bound. On the basis of these experiments, we also set $\gamma = 0$, implying that we add at most one optimality cut for each iteration.} When computing the cuts of Section \ref{subsec:MIS}, we gave the 2D-OPP check just two CPU seconds, otherwise it was given the entire remaining time limit. 

We tested two variants of the new algorithm:
\begin{itemize}
\item CHI$_{\mbox{PAOT}}$: Algorithm \ref{alg:overall} in which the initial solution is obtained by the running the metaheuristic by \citet{PAOT10};
\item CHI$_{\mbox{BKS}}$: Algorithm \ref{alg:overall} in which the initial solution is initialized to the best known solution value.
\end{itemize}
We first compare the performance of the new algorithms with that of previous methods in the literature, and then present a detailed computational analysis.

\subsection{Comparison with the Existing Literature}\label{subsec:literature}

We focus on the number of proven optimal solutions that we obtained, and compare it with those obtained by the most effective algorithms from the 2D-BPP literature. Namely, 
\begin{itemize}
\item MV98: the branch-and-bound  by \citet{MV98}, executed for 100 seconds on a Digital Alpha 533 MHz and compared with the lower bounds in \citet{MT06};
\item BM03: algorithm HBP by \citet{BM03b} in the implementation by \citet{MT06}, executed for 100 seconds on a Digital Alpha 533 MHz and compared with the lower bounds in \citet{MT06};
\item MT06: the set-covering-based heuristic by \citet{MT06}, executed for 100 seconds on a Digital Alpha 533 MHz and compared with the lower bounds in \citet{MT06};
\item CCM07: lower bound $L^2_{CCM}$ by \citet{CCM07}, compared by the authors with the best upper bound values from \citet{BM03b};
\item PS07: the dual decomposition method by \citet{PS07}, executed for one hour on an Intel Pentium IV 3.0 GHz;
\item PAOT10: the metaheuristic by \citet{PAOT10} compared  by the authors with the lower bounds in the web site of the OR group Bologna, and executed for a maximum number of iterations on a Pentium Mobile 1.5 GHz (requiring on average 15 seconds and at 70 seconds in the worst case);
\item WOZL13: the metaheuristic by \citet{WOZL13}, compared  by the authors with the lower bounds in the web site of the OR group Bologna, and executed for 2 minutes on an 
 Intel Xeon E5430 2.66 GHz Quad Core CPU.
\end{itemize}

The comparison is shown in Table \ref{tab:literature1}, grouped by class, and in Table \ref{tab:literature2}, grouped by number of items. {All numbers reported for the algorithms in the literature have been directly taken from the cited papers.}  \citet{MT06} could prove the optimality of 430 instances, including all those with 20 items. A year after, \citet{CCM07} presented their effective lower bound, which at that time achieved the highest number of optima for class 7, and \citet{PS07} proposed their exact algorithm, which remained the state-of-the-art method for the 2D-BPP until now and obtained 430 proven optima. The successive improvements obtained by the metaheuristics PAOT10 and WOZL13 increased the number of optimal solutions to 439. Both variants of our exact algorithm show good improvements with respect to the state-of-the-art. CHI$_{\mbox{PAOT}}$, which is a stand-alone algorithm that produces both lower and upper bound values, can solve 470 instances. {The use of the best known solution values allowed CHI$_{\mbox{BKS}}$ to find seven} more proven optima. It is impossible to compare precisely the CPU efforts of all methods, as some of them computed just an upper bound or used very old computers. In any case, the increase in efficiency that we obtained is remarkable. The next section better describes the contribution of each algorithmic component.

\begin{table}[htb]
\centering
\caption{Number of proven optimal solutions per class obtained by new and existing algorithms (50 instances per line, best values in bold)} \label{tab:literature1}
\scriptsize
	\begin{tabular}{r rrrrrrrrr}
	\toprule
		 & \multicolumn{7}{c}{literature} & \multicolumn{2}{c}{new}  \\
	\cmidrule(lr){2-8} \cmidrule(lr){9-10}	
	class & MV98 & BM03 & MT06 & CCM07 &  PS07 & PAOT10 & WOZL13 & CHI$_{\mbox{PAOT}}$ & CHI$_{\mbox{BKS}}$\\
	\cmidrule(lr){1-1} \cmidrule(lr){2-8} \cmidrule(lr){9-10}		
	1 & 42 & 43 & 46 & 44 & 49 & 46 & 46 & \textbf{50} & \textbf{50} \\
	2 & 43 & \textbf{50} & \textbf{50} & \textbf{50} & 48 & \textbf{50} & \textbf{50} & \textbf{50} & \textbf{50} \\
	3 & 29 & 37 & 41 & 36 & 48 & 41 & 41 & \textbf{49} & \textbf{49} \\
	4 & 39 & 45 & 45 & 44 & 44 & 45 & \textbf{49} & 46 & \textbf{49} \\
	5 & 36 & 37 & 40 & 33 & 46 & 40 & 40 & {48}  & \textbf{49} \\
	6 & 42 & 46 & 46 & 45 & 45 & 46 & \textbf{47} & 46 & \textbf{47}\\
	7 & 26 & 31 & 36 & 38 & 35 & 36 & 36 & 44  & \textbf{45}\\
	8 & 35 & 39 & 40 & 38 & 42 & 40 & 42 & \textbf{48}  & \textbf{48}\\
	9 & \textbf{50} & \textbf{50} & \textbf{50} & \textbf{50} & \textbf{50} & \textbf{50} & \textbf{50} & \textbf{50}  & \textbf{50}\\
	10 & 27 & 32 & 36 & 26 & 28 & 36 & 38 & {39}  & \textbf{40} \\
	\cmidrule(lr){1-1} \cmidrule(lr){2-8} \cmidrule(lr){9-10}	
	Total & 369 & 410 & 430 & 404 & 435 & 430 & 439 & {470}  & \textbf{477}\\	
	\bottomrule
	\end{tabular}
\end{table}

\begin{table}[htb]
\centering
\caption{Number of proven optimal solutions per number of items obtained by new and existing algorithms (100 instances per line, best values in bold)} \label{tab:literature2}
\scriptsize
	\begin{tabular}{r rrrrrrrrr}
	\toprule
		 & \multicolumn{7}{c}{literature} &  \multicolumn{2}{c}{new} \\
	\cmidrule(lr){2-8} \cmidrule(lr){9-10}	
	$n$ & MV98 & BM03 & MT06 & CCM07 &  PS07 & PAOT10 & WOZL13 & CHI$_{\mbox{PAOT}}$ & CHI$_{\mbox{BKS}}$\\
	\cmidrule(lr){1-1} \cmidrule(lr){2-8} \cmidrule(lr){9-10}		
	20 & \textbf{100} & \textbf{100} & \textbf{100} & 96 & 99 & 99 & \textbf{100} & \textbf{100}& \textbf{100} \\
	40 & 88 & 92 & 93 & 85 & 94 & 94 & 93 & \textbf{97} & \textbf{97} \\
	60 & 68 & 78 & 83 & 82 & 88 & 88 & 87 & {94}& \textbf{96} \\
	80 & 57 & 71 & 75 & 74 & 79 & 79 & 76 & {90} & \textbf{92} \\
	100 & 56 & 69 & 79 & 67 & 75 & 75 & 83 & {89} & \textbf{92} \\
	\cmidrule(lr){1-1} \cmidrule(lr){2-8} \cmidrule(lr){9-10}		
	Total & 369 & 410 & 430 & 404 & 435 & 435 & 439 & {470}& \textbf{477} \\
	\bottomrule
	\end{tabular}
\end{table}

\subsection{Detailed Computational Results}\label{subsec:detailed}

In Table \ref{tab:results}, we present the detailed computational results that we obtained with algorithm CHT$_{\mbox{BKS}}$. Each column in the table provides either an average or a sum, as next described, over the ten instances having the same class and number of items. The 'overall' line provides average or total values over the entire set of 500 instances. The columns have the following meanings:
\begin{itemize}
\item $L_c$: sum (over the 10 instances) of the continuous lower bound values, given for each instance by the sum of the item areas divided by area of the bin, without rounding up;
\item \%rmv: average (over the 10 instances) percentage of items removed from the preprocessing technique of Section \ref{subsec:preproc-pack};
\item $L'_c$: $L_c$ recomputed after the application of preprocessings \eqref{eq:min_width} and \eqref{eq:min_item_width};
\item $L_0$: sum of the $L_0$ values obtained with the methods of Section\ \ref{sec:lower_bounds};
\item $U_0$: sum of the best known upper bound values; 
\item sec$_0$: average CPU time, in seconds, required for preprocessing and lower bounds; 
\item opt$_0$: total number of proven optimal solutions with preprocessing and bounds;
\item $L$: sum of the final  lower bound values obtained;
\item $U$: sum of the final  upper bound values obtained; 
\item \#OPP: average number of calls to the 2D-OPP check procedure of Section \ref{subsec:OPP};
\item \#cuts: average number of cuts added to the model; 
\item sec$_{OPP}$: average CPU time spent in 2D-OPP checks;
\item sec:  average CPU time required by the entire execution of the algorithm;
\item opt: total number of proven optimal solutions.
\end{itemize}
The same results are also presented in Tables \ref{tab:resultsperclass} and \ref{tab:resultspern}, {where they are} grouped by, respectively, class and number of items.

From the tables, we can observe several interesting facts. The procedure to reduce the number of items performs better on small-size instances. It performs very well for some easy classes, as class 9, but much worse for other more difficult classes, as classes 4 and 6--8. With the exception of some classes, such as 4 and 6, the improvement obtained by $L'_c$ over $L_c$ is very relevant. 
Preprocessing and initial bounds are very effective in solving a large number of instances. In less than a second, on average, they can close 453 cases, including 98 of those with 20 items and the entire (although very easy) class 2. In the computation of $L_0$, $L_{BKRS}$ has always been higher than or equal to $L^2_{CCM}$, requiring however a slightly larger CPU time. 

The decomposition is able to prove 24 more optima, improving the lower bound in 23 cases and decreasing the upper bound in one case. The number of 2D-OPP checks can be very large, especially for classes 3, 7, 8 and 10. In particular, for the large instances of class 10 a few thousands checks are performed. The number of cuts added is not always related to the number of 2D-OPP checks, as some checks may have returned feasible packings and some other may have produced multiple cuts. The same consideration holds for the time spent in the 2D-OPP, which can be very high for some cases. In particular, for {an instance of class 4 with $n=80$ and another instance of} class 6 with $n=100$,  a single check consumes the entire time given to the overall algorithm, preventing it to find an optimal solution.

\begin{table}[p]
\centering
\caption{Detailed computational results of CHI$_{\mbox{BKS}}$} \label{tab:results}
\scriptsize
\renewcommand{\arraystretch}{0.85}
\setlength{\tabcolsep}{1.8mm} 
\scalebox{1.0}{
\begin{tabular}{rr rrrrrrrrrrrrrr}
\toprule
\multicolumn{2}{c}{instance} & \multicolumn{7}{c}{preprocessing and bounds} & \multicolumn{7}{c}{decomposition}\\
\cmidrule(lr){1-2} \cmidrule(lr){3-9} \cmidrule(lr){10-16}		
class & $n$ & $L_c$ & \%rmv & $L'_c$ & $L_0$ & $U_0$ & sec$_0$ & opt$_0$ & $L$ & $U$ & \#OPP & \#cuts & sec$_{OPP}$ & sec & opt\\
\cmidrule(lr){1-2} \cmidrule(lr){3-9} \cmidrule(lr){10-16}
1 & 20 & 58 & 45\% & 64 & 70 & 71 & 0.1 & 9 & 71 & 71 & 0.0 & 0.0 & 0.0 & 0.1 & 10 \\
 & 40 & 116 & 31\% & 126 & 133 & 134 & 0.1 & 9 & 134 & 134 & 0.0 & 0.0 & 0.0 & 0.2 & 10 \\
 & 60 & 179 & 17\% & 187 & 200 & 200 & 0.3 & 10 & 200 & 200 & 0.0 & 0.0 & 0.0 & 0.3 & 10 \\
 & 80 & 248 & 26\% & 261 & 275 & 275 & 0.6 & 10 & 275 & 275 & 0.0 & 0.0 & 0.0 & 0.6 & 10 \\
 & 100 & 300 & 14\% & 305 & 317 & 317 & 0.6 & 10 & 317 & 317 & 0.0 & 0.0 & 0.0 & 0.6 & 10 \\
\cmidrule(lr){1-2} \cmidrule(lr){3-9} \cmidrule(lr){10-16}
2 & 20 & 6 & 0\% & 6 & 10 & 10 & 0.1 & 10 & 10 & 10 & 0.0 & 0.0 & 0.0 & 0.1 & 10 \\
 & 40 & 13 & 0\% & 13 & 19 & 19 & 0.1 & 10 & 19 & 19 & 0.0 & 0.0 & 0.0 & 0.1 & 10 \\
 & 60 & 20 & 0\% & 20 & 25 & 25 & 0.2 & 10 & 25 & 25 & 0.0 & 0.0 & 0.0 & 0.2 & 10 \\
 & 80 & 28 & 0\% & 28 & 31 & 31 & 0.4 & 10 & 31 & 31 & 0.0 & 0.0 & 0.0 & 0.4 & 10 \\
 & 100 & 33 & 0\% & 33 & 39 & 39 & 0.5 & 10 & 39 & 39 & 0.0 & 0.0 & 0.0 & 0.5 & 10 \\
\cmidrule(lr){1-2} \cmidrule(lr){3-9} \cmidrule(lr){10-16}
3 & 20 & 38 & 44\% & 44 & 50 & 51 & 0.1 & 9 & 51 & 51 & 0.0 & 0.0 & 0.0 & 0.1 & 10 \\
 & 40 & 77 & 26\% & 85 & 92 & 94 & 0.3 & 8 & 93 & 94 & 803.0 & 1453.8 & 1.4 & 367.3 & 9 \\
 & 60 & 120 & 8\% & 121 & 138 & 139 & 0.7 & 9 & 139 & 139 & 35.2 & 17.5 & 0.0 & 64.3 & 10 \\
 & 80 & 167 & 8\% & 169 & 188 & 189 & 1.0 & 9 & 189 & 189 & 0.0 & 0.0 & 0.0 & 139.5 & 10 \\
 & 100 & 201 & 2\% & 201 & 223 & 223 & 2.4 & 10 & 223 & 223 & 0.0 & 0.0 & 0.0 & 2.4 & 10 \\
\cmidrule(lr){1-2} \cmidrule(lr){3-9} \cmidrule(lr){10-16}
4 & 20 & 6 & 0\% & 6 & 10 & 10 & 0.1 & 10 & 10 & 10 & 0.0 & 0.0 & 0.0 & 0.1 & 10 \\
 & 40 & 12 & 0\% & 12 & 19 & 19 & 0.2 & 10 & 19 & 19 & 0.0 & 0.0 & 0.0 & 0.2 & 10 \\
 & 60 & 19 & 0\% & 19 & 23 & 23 & 0.3 & 10 & 23 & 23 & 0.0 & 0.0 & 0.0 & 0.3 & 10 \\
 & 80 & 27 & 0\% & 27 & 30 & 31 & 0.5 & 9 & 30 & 31 & 57.7 & 77.1 & 359.6 & 359.8 & 9 \\
 & 100 & 32 & 0\% & 32 & 37 & 37 & 0.9 & 10 & 37 & 37 & 0.0 & 0.0 & 0.0 & 0.9 & 10 \\
\cmidrule(lr){1-2} \cmidrule(lr){3-9} \cmidrule(lr){10-16}
5 & 20 & 48 & 32\% & 56 & 65 & 65 & 0.1 & 10 & 65 & 65 & 0.0 & 0.0 & 0.0 & 0.1 & 10 \\
 & 40 & 97 & 34\% & 110 & 116 & 119 & 0.2 & 7 & 119 & 119 & 1.3 & 1.8 & 0.0 & 3.7 & 10 \\
 & 60 & 151 & 23\% & 161 & 178 & 180 & 0.5 & 8 & 180 & 180 & 0.0 & 0.0 & 0.0 & 164.4 & 10 \\
 & 80 & 210 & 18\% & 219 & 243 & 247 & 1.1 & 6 & 247 & 247 & 16.5 & 90.1 & 0.0 & 2.6 & 10 \\
 & 100 & 253 & 13\% & 259 & 280 & 281 & 1.4 & 9 & 280 & 281 & 48.0 & 178.6 & 0.0 & 360.9 & 9 \\
\cmidrule(lr){1-2} \cmidrule(lr){3-9} \cmidrule(lr){10-16}
6 & 20 & 5 & 0\% & 5 & 10 & 10 & 0.1 & 10 & 10 & 10 & 0.0 & 0.0 & 0.0 & 0.1 & 10 \\
 & 40 & 11 & 0\% & 11 & 15 & 17 & 0.3 & 8 & 15 & 17 & 0.2 & 0.0 & 720.0 & 718.1 & 8 \\
 & 60 & 17 & 0\% & 17 & 21 & 21 & 0.5 & 10 & 21 & 21 & 0.0 & 0.0 & 0.0 & 0.5 & 10 \\
 & 80 & 23 & 0\% & 23 & 30 & 30 & 0.8 & 10 & 30 & 30 & 0.0 & 0.0 & 0.0 & 0.8 & 10 \\
 & 100 & 28 & 0\% & 28 & 32 & 33 & 1.4 & 9 & 32 & 33 & 10.1 & 11.7 & 360.0 & 360.7 & 9 \\
\cmidrule(lr){1-2} \cmidrule(lr){3-9} \cmidrule(lr){10-16}
7 & 20 & 42 & 5\% & 47 & 55 & 55 & 0.0 & 10 & 55 & 55 & 0.0 & 0.0 & 0.0 & 0.0 & 10 \\
 & 40 & 91 & 2\% & 97 & 110 & 111 & 0.1 & 9 & 111 & 111 & 3.9 & 1.0 & 0.0 & 4.5 & 10 \\
 & 60 & 134 & 3\% & 139 & 157 & 158 & 0.2 & 9 & 157 & 158 & 27.3 & 117.0 & 0.0 & 359.4 & 9 \\
 & 80 & 193 & 4\% & 203 & 227 & 231 & 0.5 & 6 & 227 & 231 & 336.3 & 1377.5 & 0.0 & 1438.1 & 6 \\
 & 100 & 233 & 4\% & 241 & 271 & 271 & 1.0 & 10 & 271 & 271 & 0.0 & 0.0 & 0.0 & 1.0 & 10 \\
\cmidrule(lr){1-2} \cmidrule(lr){3-9} \cmidrule(lr){10-16}
8 & 20 & 43 & 12\% & 51 & 58 & 58 & 0.0 & 10 & 58 & 58 & 0.0 & 0.0 & 0.0 & 0.0 & 10 \\
 & 40 & 92 & 3\% & 97 & 112 & 113 & 0.1 & 9 & 113 & 113 & 2.5 & 0.4 & 0.0 & 0.7 & 10 \\
 & 60 & 136 & 3\% & 143 & 160 & 161 & 0.3 & 9 & 160 & 161 & 20.8 & 16.6 & 0.0 & 359.5 & 9 \\
 & 80 & 191 & 4\% & 199 & 223 & 224 & 0.5 & 9 & 224 & 224 & 0.0 & 0.0 & 0.0 & 0.8 & 10 \\
 & 100 & 236 & 3\% & 245 & 274 & 278 & 0.9 & 6 & 276 & 277 & 517.8 & 4228.6 & 0.1 & 416.3 & 9 \\
\cmidrule(lr){1-2} \cmidrule(lr){3-9} \cmidrule(lr){10-16}
9 & 20 & 89 & 100\% & 143 & 143 & 143 & 0.0 & 10 & 143 & 143 & 0.0 & 0.0 & 0.0 & 0.0 & 10 \\
 & 40 & 176 & 99\% & 277 & 278 & 278 & 0.1 & 10 & 278 & 278 & 0.0 & 0.0 & 0.0 & 0.1 & 10 \\
 & 60 & 272 & 96\% & 431 & 437 & 437 & 0.2 & 10 & 437 & 437 & 0.0 & 0.0 & 0.0 & 0.2 & 10 \\
 & 80 & 365 & 98\% & 575 & 577 & 577 & 0.4 & 10 & 577 & 577 & 0.0 & 0.0 & 0.0 & 0.4 & 10 \\
 & 100 & 445 & 94\% & 680 & 694 & 695 & 0.5 & 9 & 695 & 695 & 0.0 & 0.0 & 0.0 & 0.6 & 10 \\
\cmidrule(lr){1-2} \cmidrule(lr){3-9} \cmidrule(lr){10-16}
10 & 20 & 33 & 33\% & 35 & 42 & 42 & 0.0 & 10 & 42 & 42 & 0.0 & 0.0 & 0.0 & 0.0 & 10 \\
 & 40 & 63 & 17\% & 65 & 73 & 74 & 0.1 & 9 & 74 & 74 & 4.2 & 1.6 & 0.0 & 0.2 & 10 \\
 & 60 & 90 & 16\% & 92 & 98 & 100 & 0.2 & 8 & 98 & 100 & 1218.8 & 846.1 & 2.6 & 718.6 & 8 \\
 & 80 & 117 & 6\% & 118 & 124 & 128 & 0.3 & 6 & 125 & 128 & 4594.3 & 9173.6 & 15.8 & 1078.6 & 7 \\
 & 100 & 147 & 5\% & 148 & 153 & 158 & 0.5 & 5 & 153 & 158 & 5631.7 & 11598.8 & 21.1 & 1798.0 & 5 \\
\cmidrule(lr){1-2} \cmidrule(lr){3-9} \cmidrule(lr){10-16}
\multicolumn{2}{c}{overall}  & 5734 & 19\% & 6675 & 7185 & 7232 & 0.4 & 453 & 7208 & 7231 & 266.6 & 583.8 & 29.6 & 174.5 & 477 \\
\bottomrule
\end{tabular}
} 
\end{table}

Despite the relevant improvements with respect to previous literature, much needs to be done before the 2D-BPP can be considered a well-solved problem. Instances with just 40 items are still open problems more than 30 years after they were first made available on the web. Such difficulty arises from two main sources: either the instance is difficult because of the large number of combinations of tentative packings in the bins, or because of a single 2D-OPP check that cannot be solved. On the positive side, we mention that the difference between $U$ and $L$ is just one bin for all the 23 still open instances.
{We could also solve one of these reaming instances with a different configuration of CHI$_{\mbox{BKS}}$ attempted during the preliminary tests. We found indeed a 7-bin optimal solution of instance 9 of class 3 with $n$=40.}

\begin{table}[htb]
\centering
\caption{Computational results of CHI$_{\mbox{BKS}}$ per class} \label{tab:resultsperclass}
\scriptsize
\scalebox{1.0}{
\begin{tabular}{r rrrrrrrrrrrrrr}
\toprule
\multicolumn{1}{c}{instance} & \multicolumn{7}{c}{preprocessing and bounds} & \multicolumn{7}{c}{decomposition}\\
\cmidrule(lr){1-1} \cmidrule(lr){2-8} \cmidrule(lr){9-15}		
class & $L_c$ & \%rmv & $L'_c$ & $L_0$ & $U_0$ & sec$_0$ & opt$_0$ & $L$ & $U$ & \#OPP & \#cuts & sec$_{OPP}$ & sec & opt\\
\cmidrule(lr){1-1} \cmidrule(lr){2-8} \cmidrule(lr){9-15}		
1 & 901 & 26\% & 943 & 995 & 997 & 0.4 & 48 & 997 & 997 & 0.0 & 0.0 & 0.0 & 0.4 & 50 \\
2 & 100 & 0\% & 100 & 124 & 124 & 0.3 & 50 & 124 & 124 & 0.0 & 0.0 & 0.0 & 0.3 & 50 \\
3 & 603 & 17\% & 620 & 691 & 696 & 0.9 & 45 & 695 & 696 & 167.6 & 294.3 & 0.3 & 114.7 & 49 \\
4 & 96 & 0\% & 96 & 119 & 120 & 0.4 & 49 & 119 & 120 & 11.5 & 15.4 & 71.9 & 72.2 & 49 \\
5 & 759 & 24\% & 805 & 882 & 892 & 0.7 & 40 & 891 & 892 & 13.2 & 54.1 & 0.0 & 106.3 & 49 \\
6 & 84 & 0\% & 84 & 108 & 111 & 0.6 & 47 & 108 & 111 & 2.1 & 2.3 & 216.0 & 216.0 & 47 \\
7 & 695 & 3\% & 727 & 820 & 826 & 0.4 & 44 & 821 & 826 & 73.5 & 299.1 & 0.0 & 360.6 & 45 \\
8 & 699 & 5\% & 735 & 827 & 834 & 0.4 & 43 & 831 & 833 & 108.2 & 849.1 & 0.0 & 155.5 & 48 \\
9 & 1347 & 97\% & 2106 & 2129 & 2130 & 0.2 & 49 & 2130 & 2130 & 0.0 & 0.0 & 0.0 & 0.2 & 50 \\
10 & 450 & 15\% & 458 & 490 & 502 & 0.2 & 38 & 492 & 502 & 2289.8 & 4324.0 & 7.9 & 719.1 & 40 \\
\cmidrule(lr){1-1} \cmidrule(lr){2-8} \cmidrule(lr){9-15}		
\multicolumn{1}{c}{overall} & 5734 & 19\% & 6675 & 7185 & 7232 & 0.4 & 453 & 7208 & 7231 & 266.6 & 583.8 & 29.6 & 174.5 & 477 \\
\bottomrule
\end{tabular}
} 
\end{table}

\begin{table}[htb]
\centering
\caption{Computational results of CHI$_{\mbox{BKS}}$ per number of items} \label{tab:resultspern}
\scriptsize
\scalebox{1.0}{
\begin{tabular}{r rrrrrrrrrrrrrr}
\toprule
\multicolumn{1}{c}{instance} & \multicolumn{7}{c}{preprocessing and bounds} & \multicolumn{7}{c}{decomposition}\\
\cmidrule(lr){1-1} \cmidrule(lr){2-8} \cmidrule(lr){9-15}		
$n$ & $L_c$ & \%rmv & $L'_c$ & $L_0$ & $U_0$ & sec$_0$ & opt$_0$ & $L$ & $U$ & \#OPP & \#cuts & sec$_{OPP}$ & sec & opt\\
\cmidrule(lr){1-1} \cmidrule(lr){2-8} \cmidrule(lr){9-15}		
20 & 369 & 27\% & 458 & 513 & 515 & 0.1 & 98 & 515 & 515 & 0.0 & 0.0 & 0.0 & 0.1 & 100 \\
40 & 748 & 21\% & 892 & 967 & 978 & 0.2 & 89 & 975 & 978 & 81.5 & 145.9 & 72.1 & 109.5 & 97 \\
60 & 1139 & 16\% & 1330 & 1437 & 1444 & 0.3 & 93 & 1440 & 1444 & 130.2 & 99.7 & 0.3 & 166.8 & 96 \\
80 & 1569 & 16\% & 1822 & 1948 & 1963 & 0.6 & 85 & 1955 & 1963 & 500.5 & 1071.8 & 37.5 & 302.1 & 92 \\
100 & 1909 & 13\% & 2173 & 2320 & 2332 & 1.0 & 88 & 2323 & 2331 & 620.8 & 1601.8 & 38.1 & 294.2 & 92 \\
\cmidrule(lr){1-1} \cmidrule(lr){2-8} \cmidrule(lr){9-15}		
\multicolumn{1}{c}{overall} & 5734 & 19\% & 6675 & 7185 & 7232 & 0.4 & 453 & 7208 & 7231 & 266.6 & 583.8 & 29.6 & 174.5 & 477 \\
\bottomrule
\end{tabular}
} 
\end{table}

\section{Conclusion}\label{sec:conclusions}

We investigated a very interesting and challenging problem, the two-dimensional bin packing problem, which aims at packing a given set of rectangular items into the minimum set of identical rectangular bins. In order to solve it exactly, we embedded into a decomposition framework a blend of novel ideas together with some of the most successful techniques developed in the last twenty years in the area of cutting and packing. With the resulting algorithm, we managed to consistently improve the number of proven optimal solutions for the standard benchmark instances. 

Overall, we believe that nowadays, to obtain state-of-the-art results for cutting and packing problems, one should never disregard the use of preprocessing techniques, as well as initial lower and upper bounding techniques. Dual feasible functions (or conservative scales) can effectively lead to high-quality lower bounds, and some metaheuristic algorithms that have been recently proposed can produce feasible solutions that are optimal or just one-bin away from the optimum. Among the decomposition methods, both primal and dual approaches seem to work well, although in recent years the best results have been obtained by combinatorial Benders decompositions (or logic-based Benders decompositions) as the one we use in this paper.

A number of instances, some including just 40 items, still remain unsolved to proven optimality, so future research is envisaged on the problem. We also believe that it is worth investigating the (very difficult) field of three dimensional cutting and packing, including problems such as the bin packing and the strip packing. We also mention that it is interesting to determine not only the best packing, but also estimating the time that is needed to build it in practice. This research direction has not been pursued in the related literature, but can be interesting both when the packings are produced in robotized plants or by human operators.

\section*{Acknowledgments}
We thank Francisco Parre{\~n}o  for providing us with the code of \citet{PAOT10}, David Pisinger for giving us the complete results from \citet{PS07}, and Michele Monaci for sharing information on lower and upper bound values for the benchmark instances. The first author kindly acknowledges financial support from the Canadian Natural Sciences and Engineering Research Council under grant 2015-04893. The third author kindly acknowledges financial support from the University of Modena and Reggio Emilia, under grants FAR 2017 Multiscale modeling in science, industry and society and FAR 2018 Analysis and optimization of health care and pharmaceutical logistic processes. We also thank Compute Canada for providing high-performance parallel computing facilities.

\bibliographystyle{plainnat} 
\bibliography{biblio} 
\end{document}